\documentclass[12pt,a4paper]{article}
\usepackage{amsmath}
\usepackage{amssymb}
\makeatletter
    \newcommand{\rmnum}[1]{\romannumeral #1}
    \newcommand{\Rmnum}[1]{\expandafter\@slowromancap\romannumeral #1@}
  \makeatother

\newtheorem{propo}{Proposition}[section]

\newtheorem{lemma}{Lemma}[section]

\newtheorem{theo}[propo]{Theorem}

\numberwithin{equation}{section}
\begin{document}
\title{An Upper Bound for the Number of Solutions of Ternary Purely Exponential Diophantine Equations \uppercase\expandafter{\romannumeral2} }

\author{  Yongzhong Hu  and  Maohua Le }
\maketitle
\maketitle \edef \tmp {\the \catcode`@}
   \catcode`@=11
   \def \@thefnmark {}
   \@footnotetext {Supported by the National Natural Science Foundation of China(No.10971184)}
   \catcode`@=\tmp
   \let\tmp = \undefined

\begin{abstract}   Let $a,b,c$ be fixed coprime positive integers with $\min\{a,b,c\}>1$. In this paper, by analyzing the gap rule for solutions of the ternary purely exponential diophantine equation $a^x+b^y=c^z$, we prove that if $\max\{a,b,c\}\geq 10^{62}$, then the equation has at most two positive integer solutions $(x,y,z)$.
\end{abstract}



\section {Introduction}
\quad  Let $\mathbb{Z,N}$ be the sets of all integers and positive integers respectively. Let $a,b,c$ be fixed coprime positive integers with $\min\{a,b,c\}>1$. In this paper we discuss the number of solutions $(x,y,z)$ of the ternary purely exponential diophantine equation

\begin{equation}\label{1.1}
a^x+b^y=c^z, x,y,z\in\mathbb{N}.
\end{equation}

In 1933, K. Mahler \cite{mah8} used his $p-$adic analogue of the Thue-Siegel  method to prove that $(\ref{1.1})$
has only finitely many solutions $(x,y,z)$. His method is ineffective. Later, an effective result for solutions of $(\ref{1.1})$ was given by A.O. Gel'fond \cite{gel4}. Let $N(a,b,c)$ denote the number of solutions $(x,y,z)$ of $(\ref{1.1})$. As a straightforward consequence of an upper bound for the number of solutions of binary $S-$unit equations due to F. Beukers and H. P. Schlickewei \cite{beu2}, we have $N(a,b,c)\leq2^{36}$. In nearly two decades, many papers investigated the exact values of $N(a,b,c)$. The known results showed that $(\ref{1.1})$ has only a few solutions for some special cases(see the references of  \cite{hu5} and \cite{hu6} ).

Recently, Y.-Z. Hu and M.-H. Le \cite{hu5,hu6} successively proved that (\rmnum 1)   if $a,b,c$ satisfy certain divisibility conditions and $\max\{a,b,c\}$ is large enough, then $(\ref{1.1})$ has at most one solution $(x,y,z)$ with $\min\{x,y,z\}>1$. (\rmnum 2) If  $\max\{a,b,c\}>5\times 10^{27}$, then $N(a,b,c)\leq3$. R. Scott and R. Styer \cite{sco9} proved that if $2\nmid c$, then  $N(a,b,c)\leq2$. The proofs of the first two results are using the Gel'fond-Baker method with an elementary approach, and the proof of the last result is using some elementary algebraic number theory methods. In this paper, by analyzing the gap rule for solutions of  $(\ref{1.1})$ along the approach given in \cite{hu6}, we prove a general result as follows:

 \begin{theo}
If $\max\{a,b,c\}\geq 10^{62}$, then $N(a,b,c)\leq2$.
\end{theo}

 Notice that, for any positive integer $k$ with $k>1$, if $(a,b,c)=(2,2^k-1,2^k+1)$, then  $(\ref{1.1})$ has solutions $(x,y,z)=(1,1,1)$ and $(k+2,2,2)$. It implies that there exist infinitely many triples $(a,b,c)$ which make $N(a,b,c)=2$. Therefore, in general, $N(a,b,c)\leq 2$ should be the best upper bound for $N(a,b,c)$.

 \section {Preliminaries}

 \begin{lemma}\label{2l1}
 Let $t$ be a real number. If $t\geq10^{62}$, then $t>6500^6(\log t)^{18}$.
\end{lemma}

 {\bf {\it Proof.}}\ Let $F(t)=t-6500^6(\log t)^{18}$. Then we have $F'(t)=1-18\times 6500^6(\log t)^{17}/t$ and $F''(t)=18\times 6500^6(\log t)^{16}(\log t-17)/t^2$, where $F'(t)$ and $F''(t)$ are the derivative and the second derivative of $F(t)$. Since  $F'(10^{62})>0$ and  $F''(t)>0$ for $t\geq10^{62}$, we get $F'(t)>0$ for $t\geq10^{62}$. Further, since $F(10^{62})>0$, we obtain $F(t)>0$ for $t\geq10^{62}$. The  lemma is proved.$\Box$

 Let $\alpha$ be a fixed positive irrational number, and let $\alpha=[a_0,a_1,\dots]$ denote the simple continuous fraction of $\alpha$. For any nonnegative integer $i$, let $p_i/q_i$ be the $i-$th convergent of $\alpha$. By Chapter 10 of \cite{hua7}, we obtain the following two lemmas immediately.

 \begin{lemma}\label{2l2}
\begin{enumerate}
\rm \item {\it The convergents $p_i/q_i(i=0,1,\dots)$ satisfy}
$$\begin{array}{cc}
p_{-1}=1, p_0=a_0,p_{i+1}=a_{i+1}p_i+p_{i-1}, \\
q_{-1}=0, q_0=1,q_{i+1}=a_{i+1}q_i+q_{i-1},
\end{array}i\geq 0.$$
\rm \item { $p_0/q_0<p_2/q_2<\dots<p_{2i}/q_{2i}<p_{2i+2}/q_{2i+2}<\dots<\alpha$
$$<\dots<p_{2i+3}/q_{2i+3}<p_{2i+1}/q_{2i+1}<\dots<p_3/q_3<p_1/q_1, i\geq 0.$$ }
\rm \item {$1/q_i(q_{i+1}+q_i)<|\alpha-p_i/q_i|<1/q_iq_{i+1}, i\geq 0.$}
\end{enumerate}
\end{lemma}

 \begin{lemma}\label{2l3}
Let $p$ and $q$ be positive integers. If $|\alpha-p/q|<1/2q^2$, then $(p/d)/(q/d)$ is a convergent of $\alpha$, where $d=\gcd(p,q)$.
\end{lemma}

 Let $u,v,k$ be fixed positive integers such that $\min\{u,v,k\}>1$ and $\gcd(u,v)=1$.

  \begin{lemma}\label{2l4}{\rm(\cite{hu6},Lemma 4.3)}
  The equation
\begin{equation}\label{2.1}
u^l+v^m=k, l,m\in\mathbb{N}.
\end{equation}
has at most two solutions $(l,m)$.
\end{lemma}

  \begin{lemma}\label{2l5}
 Let $(l_1,m_1)$ and $(l_2,m_2)$ be two solutions of  $(\ref{2.1})$. If $l_1<l_2$, then $m_1>m_2$,
 \begin{equation}\label{2.2}
\max\{u^{l_2-l_1},v^{m_1-m_2}\}>\sqrt{k}.
\end{equation}
 and
 \begin{equation}\label{2.3}
u^{l_2-l_1}=v^{m_2}t+1,v^{m_1-m_2}=u^{l_1}t+1, t\in \mathbb N.
\end{equation}
\end{lemma}

 {\bf {\it Proof.}}\, Since
 \begin{equation}\label{2.4}
u^{l_1}+v^{m_1}=k, u^{l_2}+v^{m_2}=k,
\end{equation}
 we have

 \begin{equation}\label{2.5}
u^{l_1}\equiv-v^{m_1}\pmod k, u^{l_2}\equiv-v^{m_2}\pmod k.
\end{equation}
 If $l_1<l_2$ and $m_1\leq m_2$, then from  $(\ref{2.5})$ we get
 \begin{equation}\label{2.6}
u^{l_2-l_1}\equiv v^{m_2-m_1}\pmod k.
\end{equation}
 Since $\gcd(u,v)=1$ and $\min\{u,v\}>1$, we have $u^{l_2-l_1}\not= v^{m_2-m_1}$. Hence, by  $(\ref{2.4})$ and $(\ref{2.6})$, we get

  \begin{equation}\label{2.7}
k>\max\{u^{l_2},v^{m_2}\}>\max\{u^{l_2-l_1}, v^{m_2-m_1}\}>k,
\end{equation}
 a contradiction. Therefore, if $l_1<l_2$, then $m_1> m_2$. Moreover, by  $(\ref{2.5})$, we get $u^{l_2-l_1} v^{m_1-m_2}\equiv 1\pmod k$ and $(\ref{2.2})$.

 On the other hand, by $(\ref{2.4})$, we have

 \begin{equation}\label{2.8}
u^{l_1}(u^{l_2-l_1}-1)=v^{m_2}(v^{m_1-m_2}-1).
\end{equation}
Therefore, since $\gcd(u,v)=1$, by $(\ref{2.8})$, we get $(\ref{2.3})$. The lemma is proved.$\Box$

 \begin{lemma}\label{2l6}{\rm(\cite{ben1})}
  The equation
\begin{equation}\label{2.9}
u^l-v^m=k, l,m\in\mathbb{N}
\end{equation}
has at most two solutions $(l,m)$.
\end{lemma}

 \begin{lemma}\label{2l7}
 Let $(l_1,m_1)$ and $(l_2,m_2)$ be two solutions of $(\ref{2.9})$. If $l_1<l_2$, then $m_1<m_2$,

\begin{equation}\label{2.10}
u^{l_2-l_1}=v^{m_1}t+1,v^{m_2-m_1}=u^{l_1}t+1, t\in \mathbb N,
\end{equation}

\begin{equation}\label{2.11}
v^{m_2-m_1}>u^{l_2-l_1}>v^{m_1}
\end{equation}
and
\begin{equation}\label{2.12}
v^{m_2-m_1}>k.
\end{equation}
\end{lemma}

 {\bf {\it Proof.}}\, Since

\begin{equation}\label{2.13}
u^{l_1}-v^{m_1}=k, u^{l_2}-v^{m_2}=k,
\end{equation}
if $l_1<l_2$, then from  $(\ref{2.13})$ we get $v^{m_2}+k=u^{l_2}>u^{l_1}=v^{m_1}+k$ and $m_1<m_2$. Hence, by   $(\ref{2.13})$, we have

 \begin{equation}\label{2.14}
u^{l_1}(u^{l_2-l_1}-1)=v^{m_1}(v^{m_2-m_1}-1),
\end{equation}
whence we obtain $(\ref{2.10})$, since $\gcd(u,v)=1$. Further, by $(\ref{2.10})$ and $(\ref{2.13})$, we have

 \begin{equation}\label{2.15}
v^{m_2-m_1}-u^{l_2-l_1}=(u^{l_1}-v^{m_1})t=kt.
\end{equation}
Therefore, by $(\ref{2.10})$ and $(\ref{2.15})$, we obtain $(\ref{2.11})$ and $(\ref{2.12})$. The lemma is proved.$\Box$

Let $r,s$ be fixed coprime positive integers with $\min\{r,s\}>1$.

 \begin{lemma}\label{2l8}{\rm(\cite{car3})}
  There exist positive integers $n$ such that
\begin{equation}\label{2.16}
r^n\equiv \delta\pmod s,\delta\in\{1,-1\}.
\end{equation}
\end{lemma}
Let $n_1$ be the least value of $n$ with  $(\ref{2.16})$. Then we have $r^{n_1}\equiv \delta_1\pmod s$ and

\begin{equation}\label{2.17}
r^{n_1}=sf+\delta_1,\delta_1\in\{1,-1\}, f\in\mathbb N.
\end{equation}
A positive integer $n$ satisyies $(\ref{2.16})$ if and only if $n_1|n$. Moreover, if $n_1|n$, then $r^{n_1}-\delta_1|r^n-\delta$.

Obviously, for any fixed $r$ and $s$, the corresponding $n_1,\delta_1$ and $f$ are unique.

\begin{lemma}\label{2l9}
  Let $t$ be a positive integer such that $t>1$ and $s$ is divisble by every prime divisor of $t$. Let $n'$ be a positive integer satisfies
\begin{equation}\label{2.18}
r^{n'}\equiv\delta'\pmod {st},\delta'\in\{1,-1\}.
\end{equation}
If $s$ satisfies either $2\nmid s$ or $4|s$, then $n_1|n'$ and
\begin{equation}\label{2.19}
\frac{n'}{n_1}\equiv 0\pmod{\frac{t}{\gcd(t,f)}}.
\end{equation}
\end{lemma}

 {\bf {\it Proof.}}\,  Notice that $\gcd(r,s)=1$ and $s$ is divisble by every prime divisor of $t$. We have $\gcd(r,st)=1$. Hence, by Lemma $\ref{2l8}$, there exist positive integers $n'$ satisfy $(\ref{2.18})$. Further, since $r^{n'}\equiv\delta'\pmod s$ by $(\ref{2.18})$, we get $n_1|n'$ and
\begin{equation}\label{2.20}
n'=n_1n_2, n_2\in\mathbb N.
\end{equation}

Since either $2\nmid s$ or $4|s$, we have

\begin{equation}\label{2.21}
s>2.
\end{equation}
By  $(\ref{2.17})$,  $(\ref{2.18})$ and  $(\ref{2.20})$, we get

$$r^{n'}\equiv(r^{n_1})^{n_2}\equiv(sf+\delta_1)^{n_2}\equiv \delta_1^{n_2}+n_2\delta_1^{n_2-1}sf+$$

\begin{equation}\label{2.22}
\sum\limits_{i=2}^{n_2}\left(\begin{array}{cc}
n_2 \\
i
\end{array} \right)\delta_1^{n_2-i}(sf)^i\equiv\delta'\pmod{st}.
\end{equation}
We see from  $(\ref{2.22})$ that $\delta_1^{n_2}\equiv\delta'\pmod s$. Hence, by  $(\ref{2.21})$, we get $\delta_1^{n_2}=\delta'$, and by $(\ref{2.22})$,

\begin{equation}\label{2.23}
f\left(n_2+\sum\limits_{i=2}^{n_2}\left(\begin{array}{cc}
n_2 \\
i
\end{array} \right)(\delta_1sf)^{i-1}\right)\equiv 0\pmod t.
\end{equation}
Further, by  $(\ref{2.23})$, we obtain

\begin{equation}\label{2.24}
n_2+\sum\limits_{i=2}^{n_2}\left(\begin{array}{cc}
n_2 \\
i
\end{array} \right)(\delta_1sf)^{i-1}\equiv 0\pmod {\frac{t}{\gcd(t,f)}}.
\end{equation}

Obviously, if $t/\gcd(t,f)=1$, then $(\ref{2.19})$ holds. We just have to consider the case that $t/\gcd(t,f)>1$. Let $p$ be a prime divisor of $t/\gcd(t,f)$. Since $p|t$ and $p|s$, we see from  $(\ref{2.24})$ that $p|n_2$. Let

\begin{equation}\label{2.25}
p^\alpha\big|\big|n_2, p^\beta\big|\big|sf, p^\gamma\big|\big|\frac{t}{\gcd(t,f)}, p^{\pi_i}\big|\big|i, i\geq 2.
\end{equation}
Then, $\alpha,\beta$ and $\gamma$ are positive integers with $\beta\geq 2$ if $p=2$, $\pi_i(i\geq 2)$ are nonnegative integers satisfy

\begin{equation}\label{2.26}
\pi_i\leq\frac{\log i}{\log p}\left\{\begin{array}{cc}
\leq i-1<2(i-1)\leq\beta(i-1),  &{{\rm if}\ \ p=2,} \\
<i-1,  &{\ \ {\rm otherwise.}}
\end{array} \right.
\end{equation}
Hence, by $(\ref{2.25})$ and $(\ref{2.26})$, we have

\begin{equation}\label{2.27}
\left(\begin{array}{cc}
n_2 \\
i
\end{array} \right)(\delta_1sf)^{i-1}\equiv n_2\left(\begin{array}{cc}
n_2-1 \\
i-1
\end{array} \right)\frac{(\delta_1sf)^{i-1}}{i}\equiv0 \pmod {p^{\alpha+1}}
\end{equation}
for $i\geq2$.

By  $(\ref{2.25})$ and $(\ref{2.27})$, we get

\begin{equation}\label{2.28}
p^\alpha\big|\big|n_2+\sum\limits_{i=2}^{n_2}\left(\begin{array}{cc}
n_2 \\
i
\end{array} \right)(\delta_1sf)^{i-1}.
\end{equation}
Further, we see from $(\ref{2.24})$, $(\ref{2.25})$ and $(\ref{2.28})$ that
\begin{equation}\label{2.29}
\alpha\geq\gamma.
\end{equation}
Therefore, take $p$ through all prime divisors of $t/\gcd(t,f)$, by $(\ref{2.20})$, $(\ref{2.25})$ and $(\ref{2.29})$, we obtain $(\ref{2.19})$. The lemma is proved.$\Box$

\section {Further lemmas on solutions of $(\ref{1.1})$}

\begin{lemma}\label{3l1}{\rm(\cite{hu6},Theorem 2.1)}
 All solutions $(x,y,z)$ of $(\ref{1.1})$ satisfy $\max\{x,y,z\}<6500(\log \max \{a,b,c\})^3$.
\end{lemma}
\begin{lemma}\label{3l2}
 Let $(x,y,z)$ be a solution of $(\ref{1.1})$ with $a^{2x}<c^z$. If $b\geq3$ and $c\geq 16$, then $y/z$ is a convergent of $\log c/\log b$ with
\begin{equation}\label{3.1}
0<\frac{\log c}{\log b}-\frac{y}{z}<\frac{2}{zc^{z/2}\log b}.
\end{equation}
\end{lemma}
{\bf {\it Proof.}}\,  Since $\min\{b,c\}>1$ and $\gcd(b,c)=1$,\,$\log c/\log b$ is a positive irratrional number. Let $d=\gcd(y,z)$. Since $a^{2x}<c^z$, if $d\geq 2$, then from
$(\ref{1.1})$ we get

\begin{equation}\label{3.2}
c^{z/2}>a^x=c^z-b^y=(c^{z/d}-b^{y/d})\sum\limits_{i=0}^{d-1}c^{(d-1-i)z/d}b^{iy/d}>c^{(d-1)z/d}\geq c^{z/2},
\end{equation}
a contradiction. So we have $d=1$ and $\gcd(y,z)=1$.

Since $a^x<c^{z/2}$, we have $a^x<b^y$. Hence, by $(\ref{1.1})$, we get

\begin{equation}\label{3.3}
z\log c=\log(b^y(1+\frac{a^x}{b^y}))<y\log b+\frac{a^x}{b^y}.
\end{equation}
Since  $a^x<b^y$, by $(\ref{1.1})$, we have $c^z<2b^y$ and
\begin{equation}\label{3.4}
\frac{a^x}{b^y}<\frac{2a^x}{c^z}<\frac{2c^{z/2}}{c^z}=\frac{2}{c^{z/2}}.
\end{equation}
Hence, by $(\ref{3.3})$ and $(\ref{3.4})$, we get

\begin{equation}\label{3.5}
0<z\log c-y\log b<\frac{2}{c^{z/2}},
\end{equation}
whence we obtain $(\ref{3.1})$. On the other hand, since $b\geq3$ and $c\geq16$, we have $2/zc^{z/2}\log b<1/2z^2$. It implies that $0<\log c/\log b-y/z<1/2z^2$ by $(\ref{3.5})$. Therefore, applying Lemma $\ref{2l3}$, $y/z$ is a convergent of $\log c/\log b$ with $(\ref{3.1})$. Thus, the lemma is proved.$\Box$

Using the same method as in the proof of Lemma $\ref{3l2}$, we can obtain the following lemma immediately.

\begin{lemma}\label{3l3}
 Let $(x,y,z)$ be a solution of $(\ref{1.1})$ with $b^{2y}<c^z$. If $a\geq10^{62}$, then $x/z$ is a convergent of $\log c/\log a$ with
 \begin{equation}\label{3.6}
0<\frac{\log c}{\log a}-\frac{x}{z}<\frac{2}{zc^{z/2}\log a}.
\end{equation}
\end{lemma}

\begin{lemma}\label{3l4}
 Let $(x,y,z)$ and$(x',y',z')$ be two solutios of $(\ref{1.1})$ such that $x>x'$ and $z<z'$. If $c=\max\{a,b,c\}\geq10^{62}$, then $(y'/d)/(z'/d)$ is a convergent of $\log c/\log b$ with
\begin{equation}\label{3.7}
0<\frac{\log c}{\log b}-\frac{y'/d}{z'/d}<\frac{2}{z'ac\log b},
\end{equation}
where $d=\gcd(y',z')$.
\end{lemma}

{\bf {\it Proof.}}\, Since $x>x'$ and $z<z'$, if $a^{x'}>b^{y'}$, then we get $2a^{x'}>c^{z'}>c^z>a^x\geq a^{x'+1}\geq 2a^{x'}$, a contradiction. So we have $a^{x'}<b^{y'}$ and

\begin{equation}\label{3.8}
z'\log c=\log(b^{y'}(1+\frac{a^{x'}}{b^{y'}}))<y'\log b+\frac{a^{x'}}{b^{y'}}.
\end{equation}
Since  $2b^{y'}>c^{z'}$, we get
\begin{equation}\label{3.9}
\frac{a^{x'}}{b^{y'}}<\frac{2a^{x'}}{c^{z'}}=\frac{2}{a^{x-x'}c^{z'-z}}\cdot\frac{a^x}{c^z}<\frac{2}{ac}.
\end{equation}
Hence, by $(\ref{3.8})$ and $(\ref{3.9})$, we obtain

\begin{equation}\label{3.10}
0<\frac{\log c}{\log b}-\frac{y'}{z'}<\frac{2}{z'ac\log b}.
\end{equation}

If $|\log c/\log b-y'/z'|\geq 1/2z'^2$, then from  $(\ref{3.10})$ we get

\begin{equation}\label{3.11}
z'>\frac{1}{4}ac\log b.
\end{equation}
Since $c=\max\{a,b,c\}$, by Lemma $\ref{3l1}$, we have $z'<6500(\log c)^3$. Since $a\log b\geq\min\{2\log3,3\log2\}>2$, by $(\ref{3.11})$, we get

\begin{equation}\label{3.12}
13000(\log c)^3>c.
\end{equation}
But, since $c\geq10^{62}$, by Lemma $\ref{2l1}$, $(\ref{3.12})$ is false. Therefore, we have

\begin{equation}\label{3.13}
\left|\frac{\log c}{\log b}-\frac{y'}{z'}\right|<\frac{1}{2z'^2}.
\end{equation}
Applying Lemma  $\ref{2l3}$ to  $(\ref{3.13})$, we find from $(\ref{3.10})$ that $(y'/d)/(z'/d)$ is a convergent of $\log c/\log b$ with  $(\ref{3.7})$. Thus, the lemma is proved. $\Box$

\begin{lemma}\label{3l5}
 Let $(x,y,z)$ and$(x',y',z')$ be two solutios of $(\ref{1.1})$ such that $y>y'$ and $z\leq z'$. If $a=\max\{a,b,c\}\geq10^{62}$, then $(x'/d)/(z'/d)$ is a convergent of $\log c/\log a$ with
\begin{equation}\label{3.14}
0<\frac{\log c}{\log a}-\frac{x'/d}{z'/d}<\frac{2}{z'a\log a},
\end{equation}
where $d=\gcd(x',z')$.
\end{lemma}

{\bf {\it Proof.}}\, The proof of this lemma is similar to Lemma $\ref{3l4}$ . Since $y>y'$ and $z\leq z'$, we see from
\begin{equation}\label{3.15}
a^x+b^y=c^z,a^{x'}+b^{y'}=c^{z'}
\end{equation}
that

\begin{equation}\label{3.16}
x<x',
\end{equation}
$a^{x'}>b^{y'}$ and $2a^{x'}>c^{z'}$ . Hence, by the second equality of  $(\ref{3.15})$, we have

\begin{equation}\label{3.17}
z'\log c=\log(a^{x'}(1+\frac{b^{y'}}{a^{x'}}))<x'\log a+\frac{b^{y'}}{a^{x'}}
\end{equation}
and
\begin{equation}\label{3.18}
\frac{b^{y'}}{a^{x'}}<\frac{2b^{y'}}{c^{z'}}=\frac{2}{b^{y-y'}c^{z'-z}}\cdot\frac{b^y}{c^z}<\frac{2}{b^{y-y'}c^{z'-z}}.
\end{equation}

By  $(\ref{3.15})$ and  $(\ref{3.16})$, we have $b^y\equiv c^z\pmod {a^x}$ and $b^{y'}\equiv c^{z'}\pmod {a^{x'}}$,
whence we get
\begin{equation}\label{3.19}
b^{y-y'}c^{z'-z}\equiv 1\pmod{a^x}.
\end{equation}
Further, since $y>y'$, we have $b^{y-y'}c^{z'-z}>1$. Hence, by $(\ref{3.19})$, we get
\begin{equation}\label{3.20}
b^{y-y'}c^{z'-z}>a^x.
\end{equation}
Therefore, by  $(\ref{3.17})$,  $(\ref{3.18})$ and  $(\ref{3.20})$, we obtain $b^{y'}/a^{x'}<2/a^x$ and

\begin{equation}\label{3.21}
0<\frac{\log c}{\log a}-\frac{x'}{z'}<\frac{2}{z'a^x\log a}.
\end{equation}
Since $a^x\geq a=\max\{a,b,c\}\geq 10^{62}$, by Lemma $\ref{3l1}$, we can deduce that

\begin{equation}\label{3.22}
\frac{2}{z'a^x\log a}<\frac{1}{2z'^2}.
\end{equation}
Thus, by Lemma $\ref{2l3}$, we fond from $(\ref{3.21})$ and  $(\ref{3.22})$ that  $(x'/d)/(z'/d)$ is a convergent of $\log c/\log a$ with $(\ref{3.14})$. The lemma is proved.$\Box$

\section {The equation $A^X+\lambda B^Y=C^Z$}

\quad For any fixed triple $(a,b,c)$, put
\begin{equation}\label{4.1}
P(a,b,c)=\{(a,b,c,1),(c,a,b,-1),(c,b,a,-1)\}.
\end{equation}
 Obviously, for any element in $P(a,b,c)$, say  $(A,B,C,\lambda)$, $(\ref{1.1})$ has a solution $(x,y,z)$ is equivalent to the equation

\begin{equation}\label{4.2}
A^X+\lambda B^Y=C^Z, X,Y,Z\in\mathbb{N}
\end{equation}
has the solution

$$(X,Y,Z)=\left\{\begin{array}{cc}
(x,y,z),  &{{\rm if}\ \ (A,B,C,\lambda)=(a,b,c,1) ,} \\
(z,x,y),  &{{\rm if}\ \ (A,B,C,\lambda)=(c,a,b,-1),}\\
(z,y,x),&{{\rm if}\ \ (A,B,C,\lambda)=(c,b,a,-1).}
\end{array} \right.$$
It implies that, for any  $(A,B,C,\lambda)\in P(a,b,c)$, the numbers of solutions of $(\ref{1.1})$ and $(\ref{4.2})$ are equal. Moreover, by
Lemma $\ref{3l1}$, we have

\begin{lemma}\label{4l1}
All solutions $(X,Y,Z)$ of $(\ref{4.2})$ satisfy $\max\{X,Y,Z\}<6500(\log\max\{a,b,c\})^3$.
\end{lemma}

Here and below, we always assume that $(\ref{1.1})$ has solutions $(x,y,z)$. Then, for any $(A,B,C,\lambda)\in P(a,b,c)$, $(\ref{4.2})$ has solutions $(X,Y,Z)$.

For a fixed element $(A,B,C,\lambda)\in P(a,b,c)$, $(\ref{4.2})$ is sure to have a solution $(X_1,Y_1,Z_1)$ such that $Z_1\leq Z$, where $Z$ through all solutions $(X,Y,Z)$ of $(\ref{4.2})$ for this $(A,B,C,\lambda)$. Since $\gcd(A,C)=1$ and $\min\{A,C\}>1$, by Lemma $\ref{2l8}$, there exist positive integers $n$ such that

\begin{equation}\label{4.3}
A^n\equiv \delta\pmod{C^{Z_1}},\delta\in\{1,-1\}.
\end{equation}
Let $n_1$ be the least value of $n$ with $(\ref{4.3})$, and let
\begin{equation}\label{4.4}
A^{n_1}\equiv \delta_1\pmod{C^{Z_1}},\delta_1\in\{1,-1\}.
\end{equation}
Then we have
\begin{equation}\label{4.5}
A^{n_1}=C^{Z_1}f+\delta_1,f\in\mathbb{N}.
\end{equation}
Obviously, for any fixed triple  $(A,B,C,\lambda)\in P(a,b,c)$, the parameters $Z_1,n_1,\delta_1$ and $f$ are unique.

\begin{lemma}\label{4l2}
 $(\ref{4.2})$ has at most two solutions $(X,Y,Z)$ with the same value $Z$.
\end{lemma}

{\bf {\it Proof .}}\,\, By Lemmas $\ref{2l4}$ and $\ref{2l6}$, we obtain the lemma immediately.$\Box$

\begin{lemma}\label{4l3}{\rm (\cite{hu5}, Lemma 3.3)}
Let $(X,Y,Z)$ and $(X^\prime,Y^\prime,Z^\prime)$ be two solutions of $(\ref{4.2})$ with $Z\leq Z^\prime$. Then we have $XY^\prime-X^\prime Y\not=0$ and
$A^{|XY^\prime-X^\prime Y|}\equiv(-\lambda)^{Y+Y^\prime}\pmod {C^Z}.$
\end{lemma}

\begin{lemma}\label{4l4}.
Let $(X_1,Y_1,Z_1)$ and $(X_2,Y_2,Z_2)$ be two solutions of $(\ref{4.2})$ such that $Z_1<Z_2$ and $Z_1\leq Z$, where $Z$ through all solutions $(X,Y,Z)$ of $(\ref{4.2})$ for this
$(A,B,C,\lambda)$. If $C$ satisfies

\begin{equation}\label{4.6}
2\nmid C  \,\,\,or\,\,\,  4|C^{Z_1},
\end{equation}
then

\begin{equation}\label{4.7}
\gcd (C^{Z_2-Z_1},f)|Y_2,
\end{equation}
where $f$ is defined as in $(\ref{4.5})$.
\end{lemma}
{\bf {\it Proof .}}\,\, The proof of this lemma is similar to Lemma $\ref{2l9}$. Since $A^{X_1}+\lambda B^{Y_1}=C^{Z_1}$, $A^{X_2}+\lambda B^{Y_2}=C^{Z_2}$ and $Z_1<Z_2$, we have


$$A^{X_1Y_2}=(-\lambda)^{Y_2} B^{Y_1Y_2}+C^{Z_1}\sum\limits_{i=1}^{Y_2}\left(\begin{array}{cc}
Y_2 \\
i
\end{array} \right)(-\lambda B^{Y_1})^{Y_2-i}C^{Z_1(i-1)},$$

\begin{equation}\label{4.8}
A^{X_2Y_1}\equiv(-\lambda)^{Y_1}B^{Y_1Y_2}\pmod {C^{Z_2}}.
\end{equation}
Eliminating $B^{Y_1Y_2}$ from $(\ref{4.8})$£¬ we get

$$\lambda^\prime A^{\min\{X_1Y_2,X_2Y_1\}}\left (A^{\left |X_1Y_2-X_2Y_1\right |}-(-\lambda)^{Y_1+Y_2}\right )$$

\begin{equation}\label{4.9}
\equiv Y_2B^{Y_1(Y_2-1)}C^{Z_1}+\sum\limits_{i=2}^{Y_2}(-\lambda)^{i+1}\left(\begin{array}{cc}
Y_2 \\
i
\end{array} \right)B^{Y_1(Y_2-i)}C^{Z_1i}\pmod{C^{Z_2}},
\end{equation}
where $\lambda^\prime\in\{1,-1\}$.

 By  Lemma $\ref{4l3}$ , we have $X_1Y_2-X_2Y_1\not=0$. It implies that $\left |X_1Y_2-X_2Y_1\right|$ is a positive integer. Since $Z_1<Z_2$, using  Lemma $\ref{4l3}$ again, wer have
\begin{equation}\label{4.10}
A^{\left |X_1Y_2-X_2Y_1\right |}\equiv(-\lambda)^{Y_1+Y_2}\pmod{C^{Z_1}}.
\end{equation}
Therefore, by Lemma $\ref{2l8}$, we get from  $(\ref{4.4})$,  $(\ref{4.5})$ and  $(\ref{4.10})$ that $A^{n_1}-\delta_1\big |A^{\left |X_1Y_2-X_2Y_1\right |}-(-\lambda)^{Y_1+Y_2}$ and

\begin{equation}\label{4.11}
A^{\left |X_1Y_2-X_2Y_1\right |}-(-\lambda)^{Y_1+Y_2}=C^{Z_1}fg, g\in\mathbb N.
\end{equation}
Substitute   $(\ref{4.11})$ into   $(\ref{4.9})$, we have
$$\lambda^\prime A^{\min\{X_1Y_2,X_2Y_1\}}fg\equiv Y_2B^{Y_1(Y_2-1)}$$
\begin{equation}\label{4.12}
+\sum\limits_{i=2}^{Y_2}(-\lambda)^{i+1}\left(\begin{array}{cc}
Y_2 \\
i
\end{array} \right)B^{Y_1(Y_2-i)}C^{Z_1(i-1)}\pmod{C^{Z_2-Z_1}}.
\end{equation}

Obviously, if $\gcd(C^{Z_2-Z_1},f)=1$, then $(\ref{4.7})$ holds. We just have to consider the case that $\gcd(C^{Z_2-Z_1},f)>1$. Let $p$ be a prime divisor of $\gcd(C^{Z_2-Z_1},f)$. Since $p\big |C$ and $\gcd(B,C)=1$, we see from $(\ref{4.12})$ that $p\big | Y_2$. Let

\begin{equation}\label{4.13}
p^\alpha\big|\big|Y_2, p^\beta\big|\big|C^{Z_1}, p^\gamma\big|\big|\gcd(C^{Z_2-Z_1},f), p^{\pi_i}\big|\big|i, i\geq 2.
\end{equation}
Then, by $(\ref{4.6})$, $\alpha,\beta$ and $\gamma$ are positive integers with $\beta\geq2$ if $p=2, \pi_i(i\geq2)$ are nonnegative integers satisfy $(\ref{2.26})$. By $(\ref{2.26})$ and $(\ref{4.13})$, we have
$$\left(\begin{array}{cc}
Y_2 \\
i
\end{array} \right)B^{Y_1(Y_2-i)}C^{Z_1(i-1)}\equiv Y_2\left(\begin{array}{cc}
Y_2-1 \\
i-1
\end{array} \right)\frac{B^{Y_1(Y_2-i)}C^{Z_1(i-1)}}{i}$$
\begin{equation}\label{4.14}
\equiv0\pmod{p^{\alpha+1}}, i\geq 2.
\end{equation}
Hence, by $(\ref{4.13})$ and $(\ref{4.14})$, we get

\begin{equation}\label{4.15}
p^\alpha\big|\big|Y_2B^{Y_1(Y_2-1)}+\sum\limits_{i=2}^{Y_2}(-\lambda)^{i+1}\left(\begin{array}{cc}
Y_2 \\
i
\end{array} \right)B^{Y_1(Y_2-i)}C^{Z_1(i-1)}.
\end{equation}
Therefore, since $\gcd(C^{Z_2-Z_1},f)\big|f$ and $\gcd(C^{Z_2-Z_1},f)\big|C^{Z_2-Z_1}$, we find from $(\ref{4.12})$, $(\ref{4.13})$ and $(\ref{4.15})$ that $\alpha$ and $\gamma$ satisfies $(\ref{2.29})$. Thus, take $p$ through all prime divisors of  $\gcd(C^{Z_2-Z_1},f)$, by $(\ref{2.29})$ and $(\ref{4.13})$, we obtain $(\ref{4.7})$. The lemma is proved.$\Box$

\begin{lemma}\label{4l5}{\rm (\cite{hu6}, Lemma 4.7)}
Let $(X_j,Y_j,Z_j)(j=1,2,3)$  be three solutions of $(\ref{4.2})$ with $Z_1<Z_2\leq Z_3$. If $C=\max\{a,b,c\}$, then $\max\{a,b,c\}<5\times10^{27}.$
\end{lemma}

\begin{lemma}\label{4l6}
Let $(X_j,Y_j,Z_j)(j=1,2,3)$  be three solutions of $(\ref{4.2})$ with $Z_1<Z_2\leq Z_3$. If $C^{Z_2-Z_1}>(\max\{a,b,c\})^{1/2}$ and $C$ satisfies $(\ref{4.6})$, then $\max\{a,b,c\}<10^{62}$.
\end{lemma}
{\bf {\it Proof .}}\,\, Since $Z_2\leq Z_3$, by  Lemma $\ref{4l3}$, we have $X_2Y_3-X_3Y_2\not=0$ and

\begin{equation}\label{4.16}
A^{\left |X_2Y_3-X_3Y_2\right |}\equiv(-\lambda)^{Y_2+Y_3}\pmod{C^{Z_2}}.
\end{equation}
Further, since $Z_1<Z_2$ and $C$  satisfies $(\ref{4.6})$, by Lemma $\ref{2l9}$, we get from $(\ref{4.4})$, $(\ref{4.5})$ and $(\ref{4.16})$ that
\begin{equation}\label{4.17}
\left |X_2Y_3-X_3Y_2\right |\equiv0\pmod{\frac{C^{Z_2-Z_1}}{\gcd(C^{Z_2-Z_1},f)}},
\end{equation}
where $f$ is defined as in $(\ref{4.5})$. Recall that $X_2Y_3-X_3Y_2\not=0$. By $(\ref{4.17})$, we have
\begin{equation}\label{4.18}
\left |X_2Y_3-X_3Y_2\right |\gcd(C^{Z_2-Z_1},f)\geq C^{Z_2-Z_1}.
\end{equation}
Furthermore, by Lemma $\ref{4l4}$, we have $\gcd(C^{Z_2-Z_1},f)\leq Y_2$. Hence, we get from   $(\ref{4.18})$ that
\begin{equation}\label{4.19}
Y_2\left |X_2Y_3-X_3Y_2\right |\geq C^{Z_2-Z_1}.
\end{equation}

By Lemma $\ref{4l1}$, we have
$$Y_2\left |X_2Y_3-X_3Y_2\right |< Y_2\max\{X_2Y_3,X_3Y_2\}\leq\left(\max\{X_2,Y_2,X_3,Y_3\}\right)^3$$
\begin{equation}\label{4.20}
<6500^3\left(\log\max\{a,b,c\}\right)^9.
\end{equation}
Therefore, if $C^{Z_2-Z_1}>\left(\max\{a,b,c\}\right)^{1/2}$, then from $(\ref{4.19})$ and $(\ref{4.20})$ we get
\begin{equation}\label{4.21}
6500^6\left(\log\max\{a,b,c\}\right)^{18}>\max\{a,b,c\}.
\end{equation}
Thus, applying Lemma $\ref{2l1}$ to $(\ref{4.21})$, we obtain $\max\{a,b,c\}<10^{62}$. The lemma is proved.$\Box$

\section {Proof of Theorem 1.1 for $c=\max\{a,b,c\}$ }

\quad  By \cite{sco9}, Theorem 1.1 holds for $2\nmid c$. Therefore, we just have to consider the case that
\begin{equation}\label{5.1}
2\big|c.
\end{equation}
Since $\gcd(ab,c)=1$, by  $(\ref{5.1})$, we have

\begin{equation}\label{5.2}
2\nmid a,\,\,\,2\nmid b.
\end{equation}

In this section we will prove the theorem for the case that

\begin{equation}\label{5.3}
c=\max\{a,b,c\}\geq10^{62}.
\end{equation}
We now assume that $(\ref{1.1})$ has three solutions $(x_j,y_j,z_j)(j=1,2,3)$ with $z_1\leq z_2\leq z_3$. Then, $(\ref{4.2})$ has three solutions $(X_j,Y_j,Z_j)=(x_j,y_j,z_j)(j=1,2,3)$ for $(A,B,C,\lambda)=(a,b,c,1)$ with $Z_1\leq Z_2\leq Z_3$. By  Lemma $\ref{4l2}$, we can remove the case $z_1= z_2= z_3$. Since $C=c=\max\{a,b,c\}\geq10^{62}$, by  Lemma $\ref{4l5}$, we can remove the case $z_1<z_2\leq z_3$. So we have
\begin{equation}\label{5.4}
z_1=z_2< z_3.
\end{equation}

Since $z_1=z_2$ and

\begin{equation}\label{5.5}
a^{x_1}+b^{y_1}=a^{x_2}+b^{y_2}=c^{z_1},
\end{equation}
$(\ref{2.1})$ has two solutions $(l,m)=(x_j,y_j)(j=1,2)$ for $(u,v,k)=(a,b,c^{z_1})$. Since $(x_1,y_1)\not=(x_2,y_2)$ by $(\ref{5.5})$, we may therefore assume that
\begin{equation}\label{5.6}
x_1<x_2.
\end{equation}
Then, by  Lemma $\ref{2l5}$, we have
\begin{equation}\label{5.7}
y_1>y_2,
\end{equation}
\begin{equation}\label{5.8}
a^{x_2-x_1}=b^{y_2}t_1+1, b^{y_1-y_2}=a^{x_1}t_1+1,t_1\in\mathbb N
\end{equation}
and
\begin{equation}\label{5.9}
\max\{a^{x_2-x_1}, b^{y_1-y_2}\}>c^{z_1/2}.
\end{equation}
Accordind to the symmetry of $a$ and $b$ in  $(\ref{5.5})$, we may assume that
\begin{equation}\label{5.10}
a^{x_2-x_1}> b^{y_1-y_2}.
\end{equation}
Hence, by $(\ref{5.3})$, $(\ref{5.9})$ and $(\ref{5.10})$, we have
\begin{equation}\label{5.11}
a^{x_2-x_1}> c^{z_1/2}\geq\sqrt{c}=\left(\max\{a,b,c\}\right)^{1/2}.
\end{equation}

By $(\ref{5.6})$, if $x_3\geq x_2$, then $(\ref{4.2})$ has three solutions $(X_j,Y_j,Z_j)=(z_j,y_j,x_j)(j=1,2,3)$ for $(A,B,C,\lambda)=(c,b,a,-1)$ with
$Z_1< Z_2\leq Z_3$. Since $C^{Z_2-Z_1}=a^{x_2-x_1}> \left(\max\{a,b,c\}\right)^{1/2}$ by $(\ref{5.11})$, using Lemma $\ref{4l6}$, we get from $(\ref{5.2})$ that $\max\{a,b,c\}<10^{62}$, which contradicts $(\ref{5.3})$. Therefore, we have
 \begin{equation}\label{5.12}
x_3<x_2.
 \end{equation}

By $(\ref{5.8})$ and $(\ref{5.10})$ we get $a^{x_2-x_1}>b^{y_1-y_2}=a^{x_1}t_1+1>a^{x_1}$, and by $(\ref{5.5})$, $c^{z_1}>a^{x_2}>a^{2x_1}$. It implies that $(x,y,z)=(x_1,y_1,z_1)$ is a solution of $(\ref{1.1})$ with $a^{2x}<c^{z}$. Notice that $b\geq3$ and $c\geq16$ by $(\ref{5.2})$ and $(\ref{5.3})$. Using  Lemma $\ref{3l2}$, $y_1/z_1$ is a convergent of $\log c/\log b$ with

 \begin{equation}\label{5.13}
0<\frac{\log c}{\log b}-\frac{y_1}{z_1}<\frac{2}{z_1c^{z_1/2}\log b}.
\end{equation}

 On the other hand, by $(\ref{5.4})$ and $(\ref{5.12})$, $(x_2,y_2,z_2)$ and $(x_3,y_3,z_3)$ are two solutions of $(\ref{1.1})$ such that $x_2>x_3$ and $z_2<z_3$. Since $c=\max\{a,b,c\}\geq10^{62}$, by Lemma $\ref{3l4}$, $(y_3/d)/(z_3/d)$ is also a convergent of $\log c/\log b$ with

 \begin{equation}\label{5.14}
0<\frac{\log c}{\log b}-\frac{y_3/d}{z_3/d}<\frac{2}{z_3ac\log b},
\end{equation}
 where $d=\gcd(y_3,z_3)$.

 By $(\ref{5.4})$,$(X_1,Y_1,Z_1)=(z_1,y_1,x_1)$ and $(X_3,Y_3,Z_3)=(z_3,y_3,x_3)$ are two distinct solutions of $(\ref{4.2})$ for $(A,B,C,\lambda)=(c,b,a,-1)$. Hence, by Lemma $\ref{4l3}$, we have $z_1y_3-z_3y_1\not=0$. It implies that $y_1/z_1$ and $(y_3/d)/(z_3/d)$ are two distinct convergents of $\log c/\log b$. Therefore, by (\romannumeral2) of Lemma $\ref{2l2}$, we see from  $(\ref{5.13})$ and $(\ref{5.14})$ that

 \begin{equation}\label{5.15}
\frac{y_1}{z_1}=\frac{p_{2s}}{q_{2s}}, \frac{y_3/d}{z_3/d}=\frac{p_{2t}}{q_{2t}}, s,t\in\mathbb Z, s\not=t,\min\{s,t\}\geq0.
\end{equation}

 If $s<t$, by (\romannumeral1) and (\romannumeral3) of Lemma $\ref{2l2}$, then from  $(\ref{5.13})$ and $(\ref{5.15})$ we get
 $$z_3\geq\frac{z_3}{d}=q_{2t}\geq q_{2s+2}=a_{2s+2}q_{2s+1}+q_{2s}\geq q_{2s+1}+q_{2s}$$
 $$>\left(q_{2s}\left |\frac{\log c}{\log b}-\frac{p_{2s}}{q_{2s}}\right|\right)^{-1}=\left(z_1\left (\frac{\log c}{\log b}-\frac{y_1}{z_1}\right)\right)^{-1}$$
 \begin{equation}\label{5.16}
>\frac{1}{2}c^{z_1/2}\log b>\frac{\sqrt c}{2}.
\end{equation}
 Since $c=\max\{a,b,c\}$, by  Lemma $\ref{3l1}$, we have $z_3<6500(\log c)^3$. Therefore, by $(\ref{5.16})$, we get
 \begin{equation}\label{5.17}
13000^2\left(\log c\right)^6>c.
\end{equation}
 However, since $c\geq 10^{62}$, by Lemma $\ref{2l1}$, $(\ref{5.17})$ is false.

 Similarly, if $s>t$, then from $(\ref{5.14})$ and $(\ref{5.15})$ we get
 $$z_1=q_{2s}\geq q_{2t+2}\geq q_{2t+1}+q_{2t}>\left(q_{2t}\left |\frac{\log c}{\log b}-\frac{p_{2t}}{q_{2t}}\right|\right)^{-1}$$
 \begin{equation}\label{5.18}
=\left(\frac{z_3}{d}\left (\frac{\log c}{\log b}-\frac{y_3/d}{z_3/d}\right)\right)^{-1}>\frac{1}{2}ac\log b>c.
\end{equation}
 Further, by  Lemma $\ref{3l1}$, we have $z_1<6500(\log c)^3$. Therefore, by $(\ref{5.18})$, we get
 \begin{equation}\label{5.19}
6500\left(\log c\right)^3>c.
\end{equation}
 However, since $c\geq 10^{62}$, by Lemma $\ref{2l1}$, $(\ref{5.19})$ is false.
 Thus, we have $N(a,b,c)\leq2$ for the case $(\ref{5.3})$.

 \section {Proof of Theorem 1.1 for $c\not=\max\{a,b,c\}$ }

 \quad In this section we will prove Theorem 1.1 for the case that $c\not=\max\{a,b,c\}$. Then, by the symmetry of $a$ and $b$ in $(\ref{1.1})$, we may assume that
 \begin{equation}\label{6.1}
a=\max\{a,b,c\}\geq10^{62}.
\end{equation}
 For any solution $(x,y,z)$ of $(\ref{1.1})$, since $a^z>c^z=a^x+b^y>a^x\geq a$ by $(\ref{6.1})$, we have
 \begin{equation}\label{6.2}
z\geq2.
\end{equation}
 Hence, by $(\ref{5.1})$ and $(\ref{6.2})$, we get
 \begin{equation}\label{6.3}
4\big |c^z.
\end{equation}

 We now assume that $(\ref{1.1})$ has three solutions $(x_j,y_j,z_j)(j=1,2,3)$ with $x_1\leq x_2\leq x_3$. Then, $(\ref{4.2})$ has three solutions $(X_j,Y_j,Z_j)=(z_j,y_j,x_j)(j=1,2,3)$ for $(A,B,C,\lambda)=(c,b,a,-1)$ with $Z_1\leq Z_2\leq Z_3$. By Lemma $\ref{4l2}$, we can remove the case
 $x_1= x_2= x_3$. Since $C=a=\max\{a,b,c\}\geq10^{62}$, by Lemma $\ref{4l5}$, we an remove the case $x_1<x_2\leq x_3$. So we have
 \begin{equation}\label{6.4}
x_1=x_2<x_3.
 \end{equation}

 Since $x_1=x_2$, we have

  \begin{equation}\label{6.5}
c^{z_1}-b^{y_1}=c^{z_2}-b^{y_2}=a^{x_1}.
 \end{equation}
 It implies that $(\ref{2.9})$ has two solutions $(l,m)=(z_j,y_j)(j=1,2)$ for $(u,v,k)=(c,b,a^{x_1})$. Since $(z_1,y_1)\not=(z_2,y_2)$, we may assume that

 \begin{equation}\label{6.6}
z_1<z_2.
 \end{equation}
 Then, by Lemma $\ref{2l7}$, we get from $(\ref{6.6})$ that
 \begin{equation}\label{6.7}
y_1<y_2
\end{equation}
and
 \begin{equation}\label{6.8}
b^{y_2-y_1}=c^{z_1}t_2+1, c^{z_2-z_1}=b^{y_1}t_2+1, t_2\in\mathbb N.
\end{equation}
By the first equality of $(\ref{6.8})$, we have
 \begin{equation}\label{6.10}
b^{y_2-y_1}>c^{z_1}>a^{x_1}.
\end{equation}

 If $y_3\geq y_2$, by $(\ref{6.7})$, then $(\ref{4.2})$ has three solutions $(X_j,Y_j,Z_j)=(z_j,x_j,y_j)(j=1,2,3)$ for $(A,B,C,\lambda)=(c,a,b,-1)$ with
 $Z_1< Z_2\leq Z_3$. However, since $2\nmid b=C$ and $C^{Z_2-Z_1}=b^{y_2-y_1}>a^{x_1}\geq a=\max\{a,b,c\}\geq10^{62}$ by $(\ref{5.2})$, $(\ref{6.1})$ and $(\ref{6.10})$, using Lemma $\ref{4l6}$, it is impossible. Therefore, we obtain

 \begin{equation}\label{6.11}
y_3<y_2.
\end{equation}

If $z_3\geq z_2$, by $(\ref{6.6})$, then $(\ref{4.2})$ has three solutions $(X_j,Y_j,Z_j)=(x_j,y_j,z_j)(j=1,2,3)$ for $(A,B,C,\lambda)=(a,b,c,1)$ with
 $Z_1< Z_2\leq Z_3$. Since $a=\max\{a,b,c\}\geq10^{62}$, by Lemma $\ref{4l6}$ with $(\ref{4.3})$, we obtain

 \begin{equation}\label{6.12}
C^{Z_2-Z_1}=c^{z_2-z_1}<\left(\max\{a,b,c\}\right)^{1/2}=\sqrt a.
\end{equation}
Hence, by the second equality of  $(\ref{6.8})$ and $(\ref{6.12})$, we have
\begin{equation}\label{6.13}
b^{y_1}<b^{y_1}t_2+1= c^{z_2-z_1}<\sqrt a\leq a^{x_1/2}<c^{z_1/2}.
\end{equation}
 It implies that $(x,y,z)=(x_1,y_1,z_1)$ is a solution of $(\ref{1.1})$ with $b^{2y}<c^z$. Therefore, by Lemma $\ref{3l3}$, $x_1/z_1$ is a convergent of
 $\log c/\log a$ with
 \begin{equation}\label{6.14}
0<\frac{\log c}{\log a}-\frac{x_1}{z_1}<\frac{2}{z_1c^{z_1/2}\log a}.
\end{equation}

On the other hand, by $(\ref{6.11})$, $(\ref{1.1})$ has two solutions $(x_2,y_2,z_2)$ and $(x_3,y_3,z_3)$ such that $y_2> y_3$ and $z_2\leq z_3$. Since $a=\max\{a,b,c\}\geq10^{62}$, by Lemma $\ref{3l5}$, $(x_3/d)/(z_3/d)$ is also a convergent of
 $\log c/\log a$ with
 \begin{equation}\label{6.15}
0<\frac{\log c}{\log a}-\frac{x_3/d}{z_3/d}<\frac{2}{z_3a\log a},
\end{equation}
 where $d=\gcd(x_3,z_3)$. Further, by Lemma $\ref{4l3}$, we have $x_1z_3-x_3z_1\not=0$. It implies that $x_1/z_1$ and $(x_3/d)/(z_3/d)$ are two distinct  convergents of
 $\log c/\log a$. Hence, by (\romannumeral2) of Lemma $\ref{2l2}$, we see from  $(\ref{6.14})$ and $(\ref{6.15})$ that

 \begin{equation}\label{6.16}
\frac{x_1}{z_1}=\frac{p_{2s}}{q_{2s}}, \frac{x_3/d}{z_3/d}=\frac{p_{2t}}{q_{2t}}, s,t\in\mathbb Z, s\not=t, \min\{s,t\}\geq 0.
\end{equation}
 Since $a=\max\{a,b,c\}$, by Lemma $\ref{2l2}$ and $\ref{3l1}$, we get from $(\ref{6.14})$, $(\ref{6.15})$ and $(\ref{6.16})$ that

 $$6500\left(\log a\right)^3>\left\{\begin{array}{cc}
z_3\geq\frac{z_3}{d}=q_{2t}\geq q_{2s+2} \geq q_{2s+1}+q_{2s}\\
z_1=q_{2s}\geq q_{2t+2} \geq q_{2t+1}+q_{2t}
\end{array} \right.$$

  $$\begin{array}{cc}
>\left(q_{2s}\left |\frac{\log c}{\log a}-\frac{p_{2s}}{q_{2s}}\right|\right)^{-1}=\left(z_1\left (\frac{\log c}{\log a}-\frac{x_1}{z_1}\right)\right)^{-1}\\
>\left(q_{2t}\left |\frac{\log c}{\log a}-\frac{p_{2t}}{q_{2t}}\right|\right)^{-1}=\left(\frac{z_3}{d}\left (\frac{\log c}{\log a}-\frac{x_3/d}{z_3/d}\right)\right)^{-1}
\end{array} $$

 \begin{equation}\label{6.17}
 \left.\begin{array}{cc}
>\frac{1}{2}c^{z_1/2}\log a,  &{{\rm if}\ \ s<t} \\
>\frac{1}{2}a\log a,  &{{\rm if}\ \ s>t}
\end{array} \right\}>\frac{1}{2}\sqrt a\log a.
 \end{equation}
  But, since $a\geq10^{62}$, by Lemma $\ref{2l1}$, $(\ref{6.17})$ is false. Therefore, we obtain

 \begin{equation}\label{6.18}
 z_3<z_2.
 \end{equation}

 Finally, by the known results $(\ref{6.4})$, $(\ref{6.6})$,  $(\ref{6.7})$, $(\ref{6.11})$ and  $(\ref{6.18})$, we can complete the proof in the following four cases.

 Case \uppercase\expandafter{\romannumeral1:} $y_3\leq y_1<y_2$ and  $z_3\leq z_1<z_2$

 In this case, since  $x_3>x_1=x_2$ and

 \begin{equation}\label{6.19}
 a^{x_3}=c^{z_3}-b^{y_3}=a^{x_3-x_1}c^{z_1}-a^{x_3-x_1}b^{y_1},
 \end{equation}
we get
 \begin{equation}\label{6.20}
 c^{z_3}(a^{x_3-x_1}c^{z_1-z_3}-1)=b^{y_3}(a^{x_3-x_1}b^{y_1-y_3}-1).
 \end{equation}
 Since $\gcd(b,c)=1$, by $(\ref{6.20})$, we have
 \begin{equation}\label{6.21}
 a^{x_3-x_1}c^{z_1-z_3}=b^{y_3}t_3+1, a^{x_3-x_1}b^{y_1-y_3}=c^{z_3}t_3+1, t_3\in\mathbb N.
 \end{equation}
 Hence, we see from the second equality of $(\ref{6.21})$ that
 $a^{x_3-x_1}b^{y_1-y_3}>c^{z_3}= a^{x_3}+b^{y_3}>a^{x_3}$ and

  \begin{equation}\label{6.22}
 b^{y_1-y_3}>a^{x_1}.
 \end{equation}
 It implies that $y_3<y_1$. Therefore, $(\ref{4.2})$ has three solutions $(X_j,Y_j,Z_j)(j=1,2,3)$ for $(A,B,C,\lambda)=(c,a,b,-1)$ such that
$(X_1,Y_1,Z_1)=(z_3,x_3,y_3)$, $(X_2,Y_2,Z_2)=(z_1,x_1,y_1)$, $(X_3,Y_3,Z_3)=(z_2,x_1,y_2)$ and $Z_1< Z_2<Z_3$. However, since $2\nmid b=C$ and
 $C^{Z_2-Z_1}=b^{y_1-y_3}>a^{x_1}\geq a=\max\{a,b,c\}\geq10^{62}$, by Lemma $\ref{4l6}$, it is impossible.

  Case \uppercase\expandafter{\romannumeral2:} $y_3\leq y_1<y_2$ and  $z_1< z_3<z_2$

 Since

 \begin{equation}\label{6.23}
 c^{z_3}=a^{x_3}+b^{y_3}=c^{z_3-z_1}a^{x_1}+c^{z_3-z_1}b^{y_1},
 \end{equation}
we have
 \begin{equation}\label{6.24}
 b^{y_3}(c^{z_3-z_1}b^{y_1-y_3}-1)=a^{x_1}(a^{x_3-x_1}-c^{z_3-z_1}),
 \end{equation}
 whence we get

 \begin{equation}\label{6.25}
 c^{z_3-z_1}b^{y_1-y_3}=a^{x_1}t_4+1, a^{x_3-x_1}=b^{y_3}t_4+c^{z_3-z_1}, t_4\in\mathbb N.
 \end{equation}
 By the first equality of $(\ref{6.25})$, we have
 $c^{z_3-z_1}b^{y_1-y_3}>a^{x_1}$. It implies that

  \begin{equation}\label{6.26}
\max\{c^{z_3-z_1},b^{y_1-y_3}\}>a^{x_1/2} \geq \sqrt a.
 \end{equation}

 If $c^{z_3-z_1}>b^{y_1-y_3}$, by $(\ref{6.26})$, then $c^{z_3-z_1}>\sqrt a$. Notice that $(\ref{4.2})$ has three solutions $(X_j,Y_j,Z_j)(j=1,2,3)$ for $(A,B,C,\lambda)=(a,b,c,1)$ such that
$(X_1,Y_1,Z_1)=(x_1,y_1,z_1)$, $(X_2,Y_2,Z_2)=(x_3,y_3,z_3)$, $(X_3,Y_3,Z_3)=(x_1,y_2,z_2)$ and $Z_1< Z_2<Z_3$. However, because of $4\big |c^z=C^Z$ for any solutions $(x,y,z)$ and $(X,Y,Z)$ of $(\ref{1.1})$ and $(\ref{4.2})$ respectively, $C^{Z_2-Z_1}=c^{z_3-z_1}>\sqrt a=(\max\{a,b,c\})^{1/2}$ and $\max\{a,b,c\}\geq10^{62}$, by Lemma $\ref{4l6}$, it is impossible.

  Similarly, if $c^{z_3-z_1}<b^{y_1-y_3}$, then $y_3<y_1$ and $b^{y_1-y_3}>\sqrt a$.
  In this case, $(\ref{4.2})$ has three solutions $(X_j,Y_j,Z_j)(j=1,2,3)$ for $(A,B,C,\lambda)=(c,a,b,-1)$ such that
$(X_1,Y_1,Z_1)=(z_3,x_3,y_3)$, $(X_2,Y_2,Z_2)=(z_1,x_1,y_1)$, $(X_3,Y_3,Z_3)=(z_2,x_1,y_2)$ and $Z_1< Z_2<Z_3$. However, since $2\nmid b=C,C^{Z_2-Z_1}=b^{y_1-y_3}>\sqrt a=(\max\{a,b,c\})^{1/2}$ and $\max\{a,b,c\}\geq10^{62}$, by Lemma $\ref{4l6}$, it is impossible.

  Case \uppercase\expandafter{\romannumeral3:} $y_1<y_3<y_2$ and  $z_3\leq z_1<z_2$

 Since $x_3>x_1$, we have

 \begin{equation}\label{6.27}
 a^{x_1}+b^{y_1}=c^{z_1}\geq c^{z_3}=a^{x_3}+b^{y_3}>a^{x_1}+b^{y_1},
 \end{equation}
a contradiction.

Case \uppercase\expandafter{\romannumeral4:} $y_1<y_3<y_2$ and  $z_1< z_3<z_2$

 By  $(\ref{6.23})$, we have

 \begin{equation}\label{6.28}
 a^{x_1}(c^{z_3-z_1}-a^{x_3-x_1})=b^{y_1}(b^{y_3-y_1}-c^{z_3-z_1}).
 \end{equation}
 Since $\gcd(a,b)=1$, we get from $(\ref{6.28})$ that

\begin{equation}\label{6.29}
 c^{z_3-z_1}-a^{x_3-x_1}=b^{y_1}t_5, b^{y_3-y_1}-c^{z_3-z_1}=a^{x_1}t_5, t_5\in\mathbb N, t_5\not=0.
 \end{equation}

 If $t_5>0$, then from the second equality of $(\ref{6.29})$ we get $b^{y_3-y_1}>a^{x_1}\geq a=\max\{a,b,c\}\geq10^{62}$, In this case, $(\ref{4.2})$ has three solutions $(X_j,Y_j,Z_j)(j=1,2,3)$ for $(A,B,C,\lambda)=(c,a,b,-1)$ such that
$(X_1,Y_1,Z_1)=(z_1,x_1,y_1)$, $(X_2,Y_2,Z_2)=(z_3,x_3,y_3)$, $(X_3,Y_3,Z_3)=(z_2,x_1,y_2)$ and $Z_1< Z_2<Z_3$. However, since $2\nmid b=C,C^{Z_2-Z_1}=b^{y_3-y_1}>a=\max\{a,b,c\}\geq10^{62}$, by Lemma $\ref{4l6}$, it is impossible.

Similarly,if $t_5<0$, then from the second equality of $(\ref{6.29})$ we get $c^{z_3-z_1}>a^{x_1}\geq a=\max\{a,b,c\}\geq10^{62}$,  In this case, $(\ref{4.2})$ has three solutions $(X_j,Y_j,Z_j)(j=1,2,3)$ for $(A,B,C,\lambda)=(a,b,c,1)$ such that
$(X_1,Y_1,Z_1)=(x_1,y_1,z_1)$, $(X_2,Y_2,Z_2)=(x_3,y_3,z_3)$, $(X_3,Y_3,Z_3)=(x_1,y_2,z_2)$ and $Z_1< Z_2<Z_3$. However, since $4\big |c^z=C^Z$ for any solutions $(x,y,z)$ and $(X,Y,Z)$ of $(\ref{1.1})$ and $(\ref{4.2})$ respectively, and $C^{Z_2-Z_1}=c^{z_3-z_1}>\max\{a,b,c\}\geq10^{62}$, by Lemma $\ref{4l6}$, it is impossible.

Thus, Theorem 1.1 holds for the case $c\not=\max\{a,b,c\}$. To sum up, the theorem is proved. $\Box$

\begin{flushleft}
Yongzhong Hu\\
Department of Mathematics \\
Foshan University\\
Foshan,Guangdong 528000,China\\
E-mail:huuyz@aliyun.com
\end{flushleft}

\begin{flushleft}
Maohua Le\\
Institute of Mathematics\\
Lingnan Normal University\\
Zhanjiang,Guangdong 524048,China\\
E-mail:lemaohua2008@163.com
\end{flushleft}


\begin{thebibliography}{s2}
\bibitem{ben1} M. A. Bennett, On some exponential equation of S. S. Pillai, {\it Canad. J. Math.} 2001, {\bf 53(5):} 897-922.
\bibitem{beu2} F. Beukers aud H. P. Schlickewei, The  equation $x+y=1$ in finitely generated groups,  {\it Acta Arith.} 1996, {\bf 78(2):} 189-199.
\bibitem{car3} P. D. Carmichael, On the numerical factors of the arithmetic forms $\alpha^n\pm \beta^n$, {\it Ann. of Math.}(2), 1913, {\bf 15(1):} 30-70.
\bibitem{gel4} A. O. Gel$^\prime$fond, Sur la divisibilit${\rm\acute{e}}$ de la diff${\rm \acute{e}}$rence des puissances de deux nombres entiers par une puissance d'un id${\rm \acute{e}}$al premier, {\it Mat. Sb.} 1940, {\bf 7(1):} 7-25.
\bibitem{hu5} Y.-Z. Hu and M.-H. Le, A note on ternary purely exponential diophantine equations, {\it Acta Arith.} 2015, {\bf 171(2):} 173-182.
\bibitem{hu6} Y.-Z. Hu and M.-H. Le, A upper bound for the number of solutions of ternary purely exponential diophantine equations, {\it J. Number Theory.} 2018, {\bf 183:} 62-73.
\bibitem{hua7} L.-K. Hua, Introduction to number theorey(In Chinese), Beijing: Science Press, 1982.
\bibitem{mah8} K. Mahler, Zur Approximation algebraischer Zahlen \uppercase\expandafter{\romannumeral1}: \"{U}ber den gr\"{o}ssten Primtriler binarer formen, {\it Math.Ann.} 1933, {\bf 107:} 691-730.
\bibitem{sco9} R. Scott and R. Styer, Number of solutions to $a^x+b^y=c^z$, {\it Publ. Math. Debrecen.} 2016, {\bf 88(1-2):}  132-138.
























\end{thebibliography}
    \end{document}